\documentclass[12pt,a4paper, reqno]{amsart}
\usepackage{amsfonts,amsmath,amsthm,amssymb, amscd}
\usepackage[top=2cm, bottom=5cm, left=3cm, right=1cm]{geometry}

\newcommand{\tifrac}{\text{\raisebox{-4,5pt}[0pt][0pt]{$\widetilde{\phantom{nn}}\kern-1.1em$}} \genfrac{}{}{0pt}{}}
\newcommand{\tifracc}{\text{\raisebox{-6,5pt}[0pt][0pt]{$\widetilde{\phantom{nn}}\kern-1.1em$}} \genfrac{}{}{0pt}{}}
\newcommand{\tifraccc}{\text{\raisebox{-6,5pt}[0pt][0pt]{$\widetilde{\phantom{nnn}}\kern-1.5em$}} \genfrac{}{}{0pt}{}}
\newcommand{\tifracccc}{\text{\raisebox{-7pt}[0pt][0pt]{$\widetilde{\phantom{mmn}}\kern-2.1em$}} \genfrac{}{}{0pt}{}}

\usepackage{color}
\newtheorem{lemma}{Lemma}
\newtheorem{fact}{Fact}
\newtheorem{theorem}{Theorem}

\newtheorem{corollary}{Corollary}

\title[Generating  solutions]{Generating solutions of a linear equation   and
structure of elements of the Zelisko group}

\author[Bovdi]{V.A.~Bovdi}
\address{United Arab Emirates University, Al Ain, UAE}
\email{vbovdi@gmail.com}

\author[Shchedryk]{V.P.~Shchedryk}
\address{Pidstryhach Institute for Applied Problems of Mechanics
and Mathematics,  National Academy of Sciences of Ukraine, Lviv,  Ukraine}
\email{shchedrykv@ukr.net}


\begin{document}
\maketitle
\begin{abstract}
Solutions  of a linear equation $b=ax$ in  a homomorphic image of a commutative  B\'ezout domain  of stable range $1.5$  is developed.
It is proved that the set of solutions of a solvable linear equation contains at least one solution that divides the rest, which is
called  a generating solution. Generating solutions are pairwise  associates.   Using this result,   the structure of  elements of the
Zelisko group is investigated.
\end{abstract}

\section{ Introduction and Main results}
Let $R$ be a commutative elementary divisor ring  with $1\not=0$ (see  \cite[p.\,465]{Kaplansky}) and let $R^{n\times n}$ be the ring  of ${n\times n}$ matrices over $R$ in which $n\geq 2$. Let $U(R)$ and ${\rm GL}_n(R)$ be  groups of units of  rings $R$ and $R^{n\times n}$, respectively.  By the definition \cite[p.\,465]{Kaplansky} of elementary divisor rings,  for each  $A\in R^{n\times n}$,  there exist $P,Q\in {\rm GL}_n(R)$  (we  call them left and right transforming  matrices of the matrix $A$) such that
\begin{equation}\label{Eq:1}
PAQ={\rm diag}(\varphi _{1}, \ldots, \varphi _{k}, 0, \ldots, 0),
\end{equation}
where $\varphi _{k} \ne 0$ and  $\varphi _{i}$ is a divisor of $\varphi _{i+1}$ for  $i=1, \ldots, k-1<n$.

The diagonal matrix $\Phi:={\rm diag}(\varphi _{1}, \ldots, \varphi _{k}, 0, \ldots, 0)$ in \eqref{Eq:1} is called the {\it Smith   form}  and $\varphi _{1}, \ldots, \varphi _{k}$ are called {\it invariant factors} of the matrix $A$. Since invariant factors  in \eqref{Eq:1} are  determined uniquely up to associates, the Smith form of $A$ is defined ambiguously.

To the matrix $\Phi$  we associate   a subgroup ${\bf{G}}_{\Phi}\leq {\rm GL}_{n} (R)$ (see \cite[p.\,62]{Shch_Mon}) which is called the  {\it Zelisko group} of the matrix $\Phi$ and it  is definite  as:
\[
{\bf{G}}_{\Phi } =\{ H\in {\rm GL}_{n} (R)\mid \; \exists S\in {\rm GL}_{n} (R)\quad\text{such that}\quad  H\Phi =\Phi S\; \}.
\]
This definition  was first given  by V.~Zelisko \cite{Zel}  for the matrix over  polynomial ring $F[x]$ in which  $F$ is an algebraic closed  field of characteristic $0$.
The definition of the  Zelisko group ${\bf{G}}_{\Phi }$ over the  ring $R$  is independent of the choice of  the Smith  form $\Phi$ of  $A$ (see \eqref{Eq:1}).
Indeed, let  $\Phi_1:=\Phi\Upsilon$ in which  $\Upsilon:={\rm diag}( \varepsilon_{1}, \ldots, \varepsilon_k, 1, \ldots, 1)$ and     $ \varepsilon_i \in U(R)$. If  $H \in {\bf{G}}_{\Phi }$ then
\[
H \Phi_1=H (\Phi\Upsilon) =\Phi (S \Upsilon) =\Phi\Upsilon (\Upsilon^{-1} S\Upsilon)=\Phi_1S_1,\quad (S_1:=\Upsilon^{-1}S\Upsilon)
\]
and  ${\bf{G}}_{\Phi}\subseteq {\bf{G}}_{\Phi_1 }$.

Now, if $L \in {\bf{G}}_{\Phi_1}$,  then  $L\Phi_1=\Phi_1T$, where $T\in {\rm GL}_{n}(R)$, so
$ L\Phi=\Phi(\Upsilon T \Upsilon^{-1})$ and  ${\bf{G}}_{\Phi_1}\subseteq {\bf{G}}_{\Phi }$. Consequently, ${\bf{G}}_{\Phi_1}= {\bf{G}}_{\Phi }$.

Note that,  if $R$ is an elementary divisor domain  and   $\Phi:=\mathrm {diag}(\varphi_1, \ldots , \varphi_n)$ in \eqref{Eq:1}  is  a nonsingular  matrix (i.e.    $\det(\Phi) \ne 0$), then the group ${\bf G}_{\Phi}$
consists (see \cite[Theorem 2.6, p.\,63]{Shch_Mon}) of all invertible matrices of the following form:
\[
\left[
\begin{matrix}
  h_{11} &  h_{12} & \ldots &  h_{1.n-1} &  h_{1n} \\
   \frac{\varphi_2}{\varphi_1}h_{21} & h_{22} & \ldots & h_{2.n-1} & h_{2n} \\
  \ldots & \ldots & \ldots & \ldots & \ldots \\
   \frac{\varphi_n}{\varphi_1}h_{n1} & \frac{\varphi_n}{\varphi_2}h_{n2} & \ldots
   &  \frac{\varphi_n}{\varphi_{n-1}}h_{n.n-1} & h_{nn}
\end{matrix}
\right] \in R^{n\times n}.
\]
If for the matrix  $A$ we  fix $\Phi$ in \eqref{Eq:1}, then   the matrices $P$ and $Q$  are also defined ambiguously. As it was shown in \cite[Property 2.2, p.\,63]{Shch_Mon}, the set of left transforming matrices of $A$   coincides with the right coset  ${\bf{G}}_{\Phi} P$  of the Zelisko group ${\bf{G}}_{\Phi}$  in  ${\rm GL}_{n} (R)$. A similar property holds  for the set  of right  transforming matrices of $A$.
 Moreover,   the group ${\bf{G}}_{\Phi}$ actively used in the following.

\begin{fact}\label{Fact:1} {\cite[Theorem 4.3, p.\,138]{Shch_Mon}}
Let $R$ be a  commutative   elementary divisor domain. Let $P_{A}$ and  $P_{B}$ be  left transforming  matrices of $A, B\in R^{n\times n}$, respectively.  If $A$ and $B$ have the same Smith's  form $\Phi$, then the following conditions are equivalent:
\begin{itemize}
\item[(i)]  $A$ and $B$  are {\it right associates}, i.e. $A=BU$ for some $U\in {\rm GL}_{n} (R)$;
\item[(ii)]  $P_{B} =HP_{A} $  for some $H\in {\bf{G}}_{\Phi }$;
\item[(iii)]  ${\bf{G}}_{\Phi } P_{A} ={\bf{G}}_{\Phi } P_{B}$.
\end{itemize}
\end{fact}

We would like to note that the concept of the Zelisko group as well as its properties,   were  used by Kazimirski\u{\i} \cite{Kaz2} for the  solution of the problem of extraction of a regular divisor  of a  matrix over the  polynomial ring $F[x]$, where  $F$ is an algebraically closed  field of characteristic $0$. The properties of  the group ${\bf{G}}_{\Phi }$ in which  $\Phi\in R^{n\times n}$,   were explicitly investigated in \cite[Chapter 2.2 and  Chapter 2.6]{Shch_Mon}.

The notation $a| b$ in $R$ means that  $b=ac$ for some $c\in R$. The greatest common divisor of $a,b\in R$ is denoted  by $(a,b)$.
The  ring $R$ has {\it  stable range} $1.5$ (see \cite[p.\,961]{Shch_StablRange} and  \cite[p.\,46]{Shchedryk_3})  if for each $a, b\in R$ and  $c \in R\setminus\{0\} $   with the property  $(a,b,c) = 1$ there exists $r \in R$ such that
\[
(a + br, c) = 1.
\]
This notion  arose as  a modification of the Bass's concept of the stable range of  rings (see \cite[p.\,498]{Bass}). The examples of rings of stable range $1.5$  are Euclidean rings,  principal ideal rings,  rings of algebraic integers,
rings of integer analytic functions,  adequate  rings \cite[p.\,20]{Shch_Mon} and \cite{Bovdi_Shchedryk}. Note that the commutative rings of stable range $1.5$ coincide with rings of almost stable range $1$ \cite{Anderson_Juett, McGovern}.

If the ring $R$ has stable range $1.5$,  then some properties of  the Zelisko group  ${\bf{G}}_{\Phi }$  are  closely related to a  factorizability of the general linear group over $R$ (see \cite[Theorem 3, p.\,144]{Shch_DecompGroup} and  \cite[Chapter 2.6]{Shch_Mon}).

The subgroups of the lower- and upper- unitriangular $n\times n$-matrices of the general linear group ${\rm GL}_{n} (R)$  are  denote by  $U_{n}^{lw} (R)$ and  $U_{n}^{up} (R)$, respectively.

\begin{fact}\label{Fact:2} \cite{Shch_DecompGroup}
If $R$ is  a commutative B\'ezout domain, then the following conditions are equivalent:
\begin{itemize}
\item[(i)] $R$ has stable range  $1.5$;

\item[(ii)] ${\rm GL}_{2} (R)={\bf{G}}_{\Phi } \; U_{2}^{lw} (R)\; U_{2}^{up} (R)$ for all  $\Phi\in {\rm GL}_{2} (R)$;

\item[(iii)] ${\rm GL}_{n} (R)={\bf{G}}_{\Phi } \; U_{n}^{lw} (R)\; U_{n}^{up} (R)$ for all      $\Phi\in{\rm GL}_{n} (R) $ in which  $n\geq 2$.
\end{itemize}
\end{fact}

Note that Vaserstein and Wheland \cite{Vaserstein}  proved that if  $R$ has stable range 1, then
\[
{\rm GL}_{n} (R)=GT_{n}^{lw} (R)\; \; U_{n}^{up} (R)\; U_{n}^{lw} (R),\qquad (n\geq 2)
\]
in which $GT_{n}^{lw} (R)$ is the group of invertible lower triangular matrices.   In  the case of second order matrices over the commutative ring  $R$ the  converse  of  this statement  was proved in \cite{Nagarajan}. Moreover, it was shown  in \cite{Nagarajan} that a commutative ring $R$ is a Hermite ring of stable range 1 if and only if
\[
M_{n} (R)=T_{n}^{lw} (R)\; U_{n}^{up} (R)\; U_{n}^{lw} (R),\qquad (n\geq 2)
\]
where $T_{n}^{lw} (R)$ is the ring of lower triangular $n\times n$ matrices.   For the case of a noncommutative ring $R$,  this statement was  proved  in \cite[Theorem 1.2.2.,  p.\,12]{Chen} (see also  \cite{Chen_2}).

The notion of  rings of stable range  $1.5$ also closely related with complementability  of an unimodular row to an invertible matrix.

\begin{fact}\label{Fact:3} \cite{Shch_DecompGroup}
Let $R$  be  a  commutative B\'ezout domain. The following conditions are equivalent:
\begin{itemize}
\item[(i)] $R$ has stable range  $1.5$;

\item[(ii)]  for each relatively prime elements  $a_1, \ldots, a_n\in R$,  in which  $n\geq 3$ and  $a_1\not=0$,   there exists an invertible matrix of the following form:
\[
\left[\begin{matrix}
 {u_{n} } & {0} & {\ldots } & {0} & {0} & {u_{1} } \\
 {0} & {1} & {} & {0} & {0} & {u_{2} } \\
 {\vdots } & { } & {\ddots } & {} & {\vdots } & {\vdots } \\
  {0} & {0} & {} & {1} & {0} & {u_{n-2} } \\
   {0} & {0} & {\ldots } & {0} & {1} & {u_{n-1} } \\
    {a_{1} } & {a_{2} } & {\ldots } & {a_{n-2} } & {a_{n-1} } & {a_{n} }
     \end{matrix}\right]\in R^{n\times n}.
\]
\end{itemize}
\end{fact}

We begin our  article  by  investigating the properties of   solutions of linear equations in homomorphic images of a commutative B\'ezout domain $R$ of stable range   $1.5$.  Those  solutions of a solvable linear  equation $b=a\cdot x$ \;  ($a,b\in R$)\;   which divide  all other
are called   {\it generating solutions}  of this   equation.

 Our first result is related to generating  solutions of linear equations.
\begin{theorem}\label{Theorem:A}
Let $R$ be  a commutative B\'ezout domain (with the property $1\not=0$) of stable range $1.5$.
Let $U(R)$ be  the group of units of   $R$. For each $m\in R\setminus \{U(R), 0\}$ we denote  the factor ring $R_{m} =R/mR$.
Let $a,b\in R_m$. The following conditions hold:
\begin{itemize}
\item[(i)] each   solvable linear equation $b=ax$  in  $R_m$ has   at least one  generating solution;

\item[(ii)] each two  generating solutions of a  solvable  linear equation $b=ax$  are pairwise associates.
\end{itemize}
\end{theorem}
If we fix    an ordering relation $\leq$ on  elements of the set  $R_m$, then  the set of generating solutions of each solvable equation ${\varphi}_2={\varphi}_1\cdot {x}$ contains a minimal element which we denote  by    $\tifracc{{\varphi}_2}{{\varphi}_1}$.

Now we are able to formulate our next result.

\begin{theorem}\label{Theorem:B}
Let $R$ be  a commutative B\'ezout domain (with the property $1\not=0$) of stable range $1.5$.
Let $U(R)$ be  the group of units of   $R$. For each $m\in R\setminus \{U(R), 0\}$ we denote  the factor ring $R_{m} =R/mR$. Let   $\Phi: ={\rm diag}(\varphi _{1}, \varphi _{2},  \ldots, \varphi _{n})\in R_{m}^{n\times n}$ in which   $\varphi _{1} | \varphi _{2} |   \cdots |  \varphi _{n} \ne 0$ and $n\geq 2$.
The Zelisko group ${\bf{G}}_{\Phi }$ consists of all invertible matrices of the form:
\begin{equation} \label{Eq:2}
\left[
\begin{matrix}
{h_{11} } & {h_{12} } & {\cdots} & {h_{1,\; n-1} } & {h_{1n} } \\
\tifracc{\varphi _{2} }{\varphi _{1}} h_{21}  & {h_{22} } & {\cdots} & {h_{2,  n-1} } & {h_{2n} } \\
{\cdots} & {\cdots} & {\cdots} & {\cdots} & {\cdots} \\
\tifracc{\varphi _{n}}{\varphi _{1}} h_{n1}  & \tifracc{\varphi _{n}}{\varphi _{2}} h_{n2}  & {\cdots} &
\tifracccc{\varphi _{n} }{\varphi _{n-1}} h_{n, n-1}  & {h_{nn} } \end{matrix}
   \right]
\end{equation}
in which $h_{ij}\in R_m$ and  the element  $\tifracc{\varphi _{k}}{\varphi _{l}}\in R_m$  is the minimal  generating solution  of the linear equation ${\varphi}_{k}={\varphi}_{l}\cdot {x}$ in $R_m$ with $1\leq l<k\leq n$.
 \end{theorem}

\section{Preliminaries and Proofs}
Let $U(R)$ be   the group of units of  a commutative B\'ezout domain $R$ of stable range $1.5$.  For each $m\in R\setminus\{0, U(R)\}$ we define  the homomorphism $\overline{\bullet}: R\to R_m=R/mR$. For each $a\in R$,  we denote $\overline{a}:=\overline{\bullet}(a)\in R_m$.

We start our proof with the following.

\begin{lemma}\label{LL:1}
Let $\alpha, \beta, \sigma \in R$ such that   $\overline{\alpha}=\overline{\beta}\cdot \overline{\sigma}$.
There exist $a, b, c \in R$,  such that
\[
a=b\cdot c, \quad \overline{a}= \overline{\alpha}, \quad \overline{b}= \overline{\beta},\quad \text{and}\quad \overline{c} = \overline{\sigma}.
\]
\end{lemma}

\begin{proof}
 Set   $\overline{a}:=\alpha + mR$, $\overline{b}:=\beta + mR $ and $\overline{c}:=\sigma + mR$. Since $\beta \sigma \in \alpha +mR$,   there exists $t\in R$ such that $\beta \sigma = \alpha +mt$.
Put $a:=\alpha +mt$, $b:=\beta$ and  $c:=\sigma$.
\end{proof}

\begin{lemma}\label{LL:2}
Let  $a,b\in R$.   Elements $\overline{a}$ and $\overline{b}$ are associates  in $R_m$ if and only if $(a,m)=(b,m)$.
\end{lemma}

\begin{proof} If $\overline{a}=\overline{b}\cdot \overline{c}$, then  there are always exist  $a,b,c \in R$ such that  $a=bc$  by Lemma \ref{LL:1}.
We will use this fact freely.  Set $\mu _{a}:=(a,m)$ and   $\mu _{b}:=(b,m)$.
\smallskip

$\Leftarrow$.  Clearly,   $a=\mu _{a}  a_{1}$ and $m=\mu _{a}  m_{1}$, where   $(a_{1} ,m_{1} )=1$ and $a_{1} ,m_{1}\in R$. Thus, there exist $u,v\in R$ such that $a_{1} u+m_{1} v=1$ and for any $r\in R$ we have
\begin{equation} \label{Eq:3}
a_{1} (u+rm_{1} )+m_{1} (v-ra_{1} )=1.
\end{equation}
As $(u,m_{1} )=1$, this means that $(u,m_{1} ,\; m)=1$.  Since $R$ is a commutative B\'ezout domain of stable range $1.5$ and $m\not=0$,  there exists $r_{0} \in R$ such that $(u+r_{0} m_{1} ,\; m)=1$.  Hence $\overline{u+r_{0} m_{1}}\in U(R_{m})$. Putting  $r=r_{0} $ in \eqref{Eq:3} and multiplying by $\mu _{a}$ we get that
\begin{equation} \label{Eq:4}
a(u+r_{0} m_{1} )+m(v-r_{0} a_{1} )=\mu _{a},
\end{equation}
so $\overline{a}$ and $\overline{\mu}_{a}$ are associates  in $R_{m} $. Repeating the same reasoning, we get  $\overline{b}$ and $\overline{\mu}_{a}$ are associates  in $R_{m}$. By transitivity of the associability relation, the elements  $\overline{a}$ and $\overline{b}$ are also associates.

$\Rightarrow$.  Let $\overline{a}=\overline{b}\cdot \overline{e}$, where $\overline{e}\in U(R_{m})$. There exist $a,b\in R$, such that $a=be$ in which $(e,m)=1$ by Lemma \ref{LL:1}.  Consequently,   $(a,m)=(be,m)=(b,m)$.
\end{proof}

\begin{lemma}\label{LL:3}
Any element $\overline{a}$ in $R_m$ can be written  as  $\overline{a}=\overline{\mu} _{a}\;\overline{e}_{a}$, where  ${\mu _{a} }:=(a, m)$ is a preimage of $\overline{\mu} _{a}$, $a$ is a preimage of $\overline{a}$, and $\overline{e}_{a} \in U(R_m)$.
\end{lemma}
\begin{proof} We use notation of Lemma \ref{LL:2}.
Clearly,
$\overline{a} \; (\overline{u+r_{0} m_{1} })= \overline{(a, m)}=\overline{\mu}_{a}$
by \eqref{Eq:4}  and   $\overline{a}=\overline{\mu}_{a}\;\overline{e}_{a}$,    $\overline{e}_{a}=(\overline{u+r_{0} m_{1} })^{-1}$.
Since g.c.d. of each  element from the coset $a+mR$ with the element  $m$ is equal to $(a, m)$,  the proof is done.
\end{proof}

Note that the presentation  of $\overline{a}\in R_m$ in Lemma \ref{LL:3}  in the form  $\overline{a}=\overline{\mu} _{a}\;\overline{e}_{a}$ is ambiguous.
\smallskip

\noindent {\bf Example 1}. The element $\overline{4}\in \mathbb{Z}_{6}$ can be write as   $\overline{4}=\overline{2}\cdot\overline{2}=\overline{2}\cdot\overline{5}$, where ${2}=({4}, {6})$ and  $\overline{5} \in U(\mathbb{ Z}_{6})$, but $\overline{2}\not \in U(\mathbb{ Z}_{6})$.
Furthermore,  $ \mathbb{Z}_{36}\ni \overline{8}= \overline{4}\cdot\overline{11}=\overline{4}\cdot\overline{29}$ in which  $\overline{11}, \overline{29} \in U(\mathbb{ Z}_{36})$.

\begin{lemma}\label{LL:4}
If  $\overline{a}, \overline{b}\in R_m$ are  multiples of each others,  then $\overline{a}$ and  $\overline{b}$ are  associates in $R_m$.
\end{lemma}

\begin{proof}
If $\overline{a}=\overline{b}\:\overline{c}$, then  there  exist  $a,b,c \in R$ such that  $a=bc$  by Lemma \ref{LL:1}. Set $\mu _{a}:=(a,m)$ and   $\mu _{b}:=(b,m)$. Clearly, $\mu _{a}=(bc,m)$, so  $\mu _{b} | \mu _{a}$.

Similarly,  from  $\overline{b}=\overline{a}\:\overline{d}$ follows that  $\mu _{a} | \mu _{b}$.  Since $R$ is a domain,
$\mu _{a} = \mu _{b} e$, where   $e \in U(R)$. Furthermore,
$\overline{a}=\overline{\mu}_{a}\cdot\overline{e}_{a}$ and   $\overline{b}=\overline{\mu} _{b}\cdot \overline{e}_{b}$, where $\overline{e}_{a}, \overline{e}_{b} \in U(R_m)$ by Lemma \ref{LL:3}, so
\[
\overline{a}=\overline{\mu}_{a}\;\overline{e}_{a}=(\overline{\mu}_{b}\cdot \overline{e} )\; \overline{e}_{a}
=\overline{\mu}_{b}\;\overline{e}_{b} (\overline{e}_{b})^{-1} \overline{e} \; \overline{e}_{a}
=\overline{b}\;\overline{\xi},
\]
where $\overline{\xi}= (\overline{e}_{b})^{-1} \overline{e} \; \overline{e}_{a} \in U(R_m)$.
\end{proof}

Let  $\overline{c}\in R_m$. The annihilator of $\overline{c}$  in $R_m$ is denoted  by ${\rm Ann}(\overline{c})$.

\begin{lemma}\label{LL:5}
If $b \in R$,  then ${\rm Ann}(\overline{b})=\overline{\alpha}_{b}R_{m}$, where $\alpha _{b}:=\frac{m}{\mu _{b}}\in R$ and $\mu _{b}:=(b,m)$.
\end{lemma}

\begin{proof} The ideal ${\rm Ann}(\overline{b})$ consists of the images of   $s\in R$ for which $bs=mp$, so
\[
\textstyle \frac{b}{\mu_b}s=\frac{m}{\mu_b} p\qquad \text{and}\qquad \frac{m}{\mu _b} | \frac{b}{\mu _b} s.
\]
We conclude that  $\frac{m}{\mu _b} | s$ and $s=\frac{m}{\mu _b} s'$ for some $ s'\in R$ because $\left(\frac{m}{\mu _b} ,\frac{b}{\mu _b} \right)=1$. Consequently,  $\overline{c}=\overline{\alpha}_{b}\cdot \overline{s'}\in \overline{\alpha}_{b}R_{m} $ and  ${\rm Ann}(\overline{b})\subseteq \overline{\alpha}_{b}R_{m} $.

Since $b\left(\frac{m}{\mu _b} t\right)=m\left(\frac{b}{\mu _b} t\right)$ for each $t\in R$,  the image of $\frac{m}{(b,m)} t$ lies in ${\rm Ann}(\overline{b})$. This yields  that $\overline{\alpha}_{b}R_{m} \subseteq {\rm Ann}(\overline{b})$, so  ${\rm Ann}(\overline{b})=\overline{\alpha}_{b}R_{m} $.
\end{proof}

\begin{proof}[\underline{Proof of Theorem \ref{Theorem:A}}]
(i) Let $\overline{a}, \overline{b}, \overline{c}\in R_m\setminus\{0\}$ such that $\overline{a}=\overline{b}\overline{c}$ and   $a=bc$  (see Lemma \ref{LL:1}). Thus  $\overline{a}=\overline{\mu}_{a}\cdot \overline{e}_{a}$ and $ \overline{b}=\overline{\mu}_{b}\cdot \overline{e}_{b}$, where $\mu _{a}:=(a,m)$,  $\mu _{b}:=(b,m)$ and  $\overline{e}_{a},\overline{e}_{b}\in U(R_{m} )$ by Lemma \ref{LL:3}. This yields
\[
\textstyle
\frac{\mu_{a}}{\mu _{b} } =\frac{(a,m)}{(b,m)} =\frac{(bc,m)}{(b,m)} =\left(\frac{bc}{(b,m)} ,\; \frac{m}{(b,m)} \right)=\left(\frac{b}{(b,m)} c,\; \frac{m}{(b,m)} \right)=\left(c,\; \frac{m}{(b,m)} \right),
\]
so   $\mu _{a} =\mu _{b} \sigma$, where   $\sigma:=(c,\; \frac{m}{(b,m)} )$.   It is easy to check that
\[
\overline{a}=\overline{\mu}_{a}\cdot \overline{e}_{a}=(\overline{\mu}_{b}\cdot \overline{e}_{b})(\overline{\sigma }\cdot (\overline{e}_{b})^{-1}\cdot  \overline{e}_{a})=\overline{b}\cdot \overline{c_{1} },
\]
where $\overline{c_{1} }=\overline{\sigma }\cdot (\overline{e}_{b})^{-1}\cdot  \overline{e}_{a}$.  It follows that the set of solutions of the equation
 $\overline{a}=\overline{b}\cdot \overline{x}$ is
\[
\overline{c}_{1} +{\rm Ann}(\overline{b})=
\overline{\sigma }\cdot (\overline{e}_{b})^{-1}\cdot  \overline{e}_{a}  +{\rm Ann}(\overline{b}).
\]
Thus ${\rm Ann}(\overline{b})=\overline{\alpha}_{b}R_{m}$, where $\alpha_b:=\frac{m}{\mu _{b} }$ by  Lemma \ref{LL:5}. As $\sigma =\frac{\mu_{a}}{\mu _{b} }$, in which  $\mu _{a} | m$, we get  $\sigma | \frac{m}{\mu _{b} } R$. Therefore $\overline{\sigma }| {\rm {\rm Ann}}(\overline{b})$ and $\overline{c}_{1} =\overline{\sigma }\cdot (\overline{e}_{b})^{-1}\cdot  \overline{e}_{a}   | {\rm Ann}(\overline{b})$, so  $\overline{c}_{1} | \big(\overline{c}_{1} +{\rm Ann}(\overline{b})\big)$.  Hence $\overline{c}_{1} $ is a solution of $\overline{a}=\overline{b}\cdot \overline{x}$, and the divisor of all of the other solutions of this equation.
\end{proof}

Note that,  the  solution  of a solvable linear equation $\overline{a}=\overline{b}\cdot \overline{x}$ in $R_m$  which divides all other solutions
is called   {\it generating solution}  of this equation.

\begin{proof}[\underline{Proof of Theorem \ref{Theorem:A}(ii)}] Let $f,g$ be generating solutions of a linear  equation $b=ax$. It follows that  $f|g$ and $g|f$. From  Lemma  \ref{LL:4}, we conclude that   $f,g$ are associates.\end{proof}

\smallskip

\noindent
{\bf Example 2.} Let $R_{m} =\mathbb{ Z}_{36}$.  The set $\overline{6}+{\rm Ann}(\overline{4})=\{ \overline{6},\overline{15},\overline{24},\overline{33}\} $
consists  of all  solutions of the solvable equation $\overline{4}\overline{x}=\overline{24}$, where
${\rm Ann}(\overline{4})=\{ \overline{0},\overline{9},\overline{18},\overline{27}\}$. The elements $\overline{15}$ and  $\overline{33}$
are generating solutions of our equation. These elements  divide all elements from $\overline{6}+{\rm Ann}(\overline{4})$ and are pairwise  associates, because  $\overline{33}=\overline{15}\cdot \overline{31}=\overline{15}\cdot \overline{7}$,  where $\overline{7}, \overline{31}\in U(\mathbb{ Z}_{36})$.

Note that the generating solutions of a linear equation can be characterized as the g.c.d. of all  solutions of this equation. However, in general, the  g.c.d. of two  solutions is not  a solution.  Indeed, let $R_m=\mathbb{ Z}_{72}$. The numbers $\overline{2}, \overline{ 20  }, \overline{38  }, \overline{ 56}$ are all solutions of the equation $\overline{4}\overline{ x }=\overline{ 8 }$.
Obviously,       $\overline{2}=(\overline{2}, \overline{ 20}, \overline{38}, \overline{ 56})$ are  solutions of our  equation,  but $\overline{4} =(\overline{ 20}, \overline{ 56})$ is not a solution.

\smallskip
Recall that,  if we fix    an ordering relation $\leq$ on  elements of the set  $R_m$, then  the set of generating solutions of each solvable equation $\overline{\varphi}_2=\overline{\varphi}_1\cdot \overline{x}$ contains a minimal element which we denote  by    $\tifracc{\overline{\varphi}_2}{\overline{\varphi}_1}$.
\begin{corollary}\label{Cor:1}
If  $\overline{\varphi}_1, \overline{\varphi}_2\in R_m$ such that ${\overline{\varphi}_1}\; | \; {\overline{\varphi}_2}\not=0$, then
  $\overline{\varphi}_2=\tifracc{\overline{\varphi}_2}{\overline{\varphi}_1}\cdot \overline{\varphi}_1$.
\end{corollary}
\begin{proof}
Clearly $\tifracc{\overline{\varphi}_2}{\overline{\varphi}_1}$ is a solution of the equation
$\overline{\varphi}_2=\overline{\varphi}_1\cdot \overline{x}$, so  $\overline{\varphi}_2=\tifracc{\overline{\varphi}_2}{\overline{\varphi}_1}\cdot \overline{\varphi}_1$.
\end{proof}
\smallskip

Let ${\overline{\varphi}_1}, {\overline{\varphi}_2},  \ldots , {\overline{\varphi}_n} \in R_{m}$ such that  ${\overline{\varphi}_1}\; | \; {\overline{\varphi}_2}\; | \; \cdots \; | \; {\overline{\varphi}_n}\not=0$.   Define the following sets:
\begin{equation}\label{Eq:5}
\overline{ M}_{ij}:=\{ \overline{x}\in R_m \mid {\overline{\varphi}_i}= {\overline{\varphi}_j}\cdot\overline{x}\}, \qquad\qquad (1\leq j<i\leq n).
\end{equation}
\smallskip

\begin{lemma}\label{LL:6}
For all $\overline{\mu}_{ij} \in \overline{M}_{ij}$ ($1\leq j<i\leq n$) (see \eqref{Eq:5})   the following holds:
\[
\overline{x}_{p,k}:=\overline{ \mu}_{p,p-1}\cdot \overline{\mu}_{p+1,p} \cdots  \overline{ \mu}_{p+k,p+k-1}\in  \overline{ M}_{p+k, p-1}, \qquad  (1<p\leq n, \;\; 0\leq k \leq n-p).
\]
\end{lemma}

\begin{proof} Indeed, using the definition of $\overline{ M}_{ij}$, we have
\[
\begin{split}
\overline{{\varphi}}_{p-1} \cdot  \overline{x}_{p,k} &=
(\overline{{\varphi}}_{p-1} \cdot \overline{ \mu}_{p,p-1})\; \overline{ \mu}_{p+1,p} \cdots \overline{ \mu}_{p+k,p+k-1}\\
&=({\overline{{\varphi}_p}}\cdot   \overline{ \mu}_{p+1,p})\; \overline{ \mu}_{p+2,p+1} \ldots \overline{ \mu}_{p+k,p+k-1}\\
& = \cdots= \\
&=\overline{{\varphi}}_{p+k-1} \cdot \overline{ \mu}_{p+k,p+k-1}\\
&={\overline{\varphi}_{p+k}}.
\end{split}
\]
Consequently, $\overline{x}_{p,k}\in \overline{ M}_{p+k, p-1}.$
\end{proof}

\begin{lemma}\label{LL:7}
Let ${\overline{\varphi}_1}, {\overline{\varphi}_2} , {\overline{\varphi}_3} \in R_{m}$. If  ${\overline{\varphi}_1}\; | \; {\overline{\varphi}_2}\; | \; {\overline{\varphi}_3}\not=0$, then
\[
\textstyle \textstyle \tifracc{\overline{\varphi}_2}{\overline{\varphi}_1}\cdot
\tifracc{\overline{\varphi}_3}{\overline{\varphi}_2} = \tifracc{\overline{\varphi}_3}{\overline{\varphi}_1} \cdot \overline{e }\in \overline{M}_{31}, \qquad \qquad (\overline{e }\in U(R_{m})).
\]
Moreover, $\tifracc{\overline{\varphi}_3}{\overline{\varphi}_1} \cdot \overline{e }$ is a generating solution of \; $\overline{\varphi}_3=\overline{\varphi}_1\overline{x}$.
\end{lemma}

\begin{proof} There exist  ${\varphi}_1, \varphi_2, {\varphi}_3\in R$ which are  preimages of ${\overline{\varphi}_1}, {\overline{\varphi}_2} , {\overline{\varphi}_3}$, such that
${\varphi_1} \; |   {\varphi_2}    \; |  {\varphi_3}\not=0$ by Lemma \ref{LL:1}.  Clearly   $\overline{\varphi}_i=\overline{\mu}_{\varphi_i}\;\overline{e}_{\varphi_i}$, in which   $\overline{e}_{\varphi_i} \in U(R_{m} )$ and   $\mu_{\varphi_i}:=(\varphi_i, m)$ by Lemma \ref{LL:3}.
Using the same argument,  as  in the proof of Theorem \ref{Theorem:A}(i), we get  that each  generating solution  of  the linear equation $\overline{\varphi}_i=\overline{\varphi}_j\; {\overline{x }}$   has the following form
\[
\textstyle
\overline{\psi}_{ij}:=
\overline{\left(
\frac{\mu_{\varphi_i}}{\mu_{\varphi_j}}
\right)}
\cdot \overline{e}_i \cdot (\overline{e}_j )^{-1},
\]
in which $\overline{e}_i, \overline{e}_j\in U(R_m)$ and $1\leq j<i\leq n$.
Since ${\varphi_1} \; |   {\varphi_2}    \; | {\varphi_3}\not=0$,
\[
\textstyle \frac{{\mu_{\varphi_2}}}{{\mu_{\varphi_1}}} \cdot \frac{{\mu_{\varphi_3}}}{{\mu_{\varphi_2}}}=
\frac{{\mu_{\varphi_3}}}{{\mu_{\varphi_1}}}\quad  \text {and}\quad
\textstyle
\overline{\left(\frac{{\mu_{\varphi_2}}}{{\mu_{\varphi_1}}}\right)}   \cdot  \overline{\left(\frac{{\mu_{\varphi_3}}}{{\mu_{\varphi_2}}}\right)}=
 \overline{\left(
\frac{{\mu_{\varphi_3}}}{{\mu_{\varphi_1}}}
\right)}.
\]
This yields that
\begin{equation}\label{Eq:6}
\overline{\psi}_{21} \;\overline{\psi}_{32}=\overline{\psi}_{31}.
\end{equation}
According to Theorem \ref{Theorem:A} (ii),
\[
\textstyle
\tifracc{\overline{\varphi}_2}{\overline{\varphi}_1}=\overline{\psi}_{21}\cdot \overline{\varepsilon}_{21}, \qquad  \tifracc{\overline{\varphi}_3}{\overline{\varphi}_2}=\overline{\psi}_{32} \cdot \overline{\varepsilon}_{32}, \qquad
\tifracc{\overline{\varphi}_3}{\overline{\varphi}_1}=\overline{\psi}_{31} \cdot \overline{\varepsilon}_{31},\qquad (\overline{\varepsilon}_{ij}\in U(R_{m}))
\]
and $\textstyle \tifracc{\overline{\varphi}_2}{\overline{\varphi}_1}  \cdot \tifracc{\overline{\varphi}_3}{\overline{\varphi}_2} = \tifracc{\overline{\varphi}_3}{\overline{\varphi}_1} \cdot \overline{e }$ where $
\overline{e }:=(\overline{\varepsilon}_{31})^{-1}\cdot \overline{\varepsilon}_{21}\cdot \overline{\varepsilon}_{32}
\in U(R_{m})$. Since $\tifracc{\overline{\varphi}_2}{\overline{\varphi}_1}  \cdot \tifracc{\overline{\varphi}_3}{\overline{\varphi}_2}\in\overline{ M}_{31}$  by Lemma \ref{LL:6},
$\tifracc{\overline{\varphi}_3}{\overline{\varphi}_1} \cdot \overline{e } \in\overline{ M}_{31}$.
The element $ \tifracc{\overline{\varphi}_3}{\overline{\varphi}_1}$ is a generating solution of $\overline{\varphi}_3=\overline{\varphi}_1\overline{x}$.
According to Theorem \ref{Theorem:A}(ii), the element $\tifracc{\overline{\varphi}_3}{\overline{\varphi}_1} \cdot \overline{e }$ is also  a generating solution of the same  equation. \end{proof}

Proving the previous lemma, we obtained \eqref{Eq:6}, which can be formulated as the following independent result.

\begin{corollary}\label{Cor:2}
Each  set $\overline{ M}_{ij}$  \rm{(see \eqref{Eq:5})} contains a   generating solution $\overline{ \psi}_{ij}$  such that
\[
\overline{ \psi}_{i,i-1}\cdot \overline{ \psi}_{i+1,i}=\overline{ \psi}_{i+1,i-1},\qquad\qquad (2\leq i <n-1).
\]
\end{corollary}
\hfill $\Box$
\smallskip

\noindent
{\bf Example 3.} Note that, a product of solutions  of  two equations $\overline{\varphi}_2=\overline{\varphi}_1\;\overline{x}$ and
$\overline{\varphi}_3=\overline{\varphi}_2\;\overline{x}$, in which at least one  factor  is a non generated solution, in general  is not a generating solution of $\overline{\varphi}_3=\overline{\varphi}_1\;\overline{x}$.

Indeed, let $R_m=\mathbb{ Z}_{72}$  and let  $\{\overline{\varphi}_1, \overline{\varphi}_2, \overline{\varphi}_3\}=\{\overline{4}, \overline{8}, \overline{24}\}$.
\[
\begin{array}{| c| c| c| }
  \hline
  \mbox{Equation} & \mbox{Solutions} &  \mbox{Generating solutions} \\
  \hline
  \overline{4}\overline{ x }=\overline{ 8 } &  \overline{2  }, \overline{ 20  }, \overline{38  }, \overline{ 56  }  &  \overline{\bf  2 }, \overline{ 38 } \\
  \hline
   \overline{8}\overline{ x }=\overline{ 24 } &   \overline{3  }, \overline{ \bf 12  }, \overline{ 21  }, \overline{30  },  \overline{39 }, \overline{48  }, \overline{ 57  }, \overline{ 66  } &  \overline{3  }, \overline{21 },  \overline{39 },  \overline{57 }  \\
   \hline
   \overline{4}\overline{ x }=\overline{24 } &  \overline{6  }, \overline{24  }, \overline{42  }, \overline{ 60  } &  \overline{ 6 }, \overline{ 42 } \\
  \hline
\end{array}
\]
However,   $\overline{2 } \cdot \overline{12}=\overline{24}$ is not a generating solution of the equation
$\overline{24}={\overline{4}}\cdot \overline{x}$.

\begin{lemma}\label{LL:8}
Let $\overline{a}, \overline{b}\in R_m$.
If ${\overline{b}}\; | \;{\overline{a}}$ then \;
$\tifrac{\overline{a}}{\overline{b}}\;  =\tifracc{\overline{\alpha} _{b}}{\overline{\alpha} _{a}} \cdot \overline{e}$ in which
\[
\textstyle
\overline{\alpha}_{a}:=\overline{\left(\frac{m}{(a, m)}\right)},
\quad
\overline{\alpha}_{b}:=\overline{\left(\frac{m}{(b, m)}\right)},
\quad \text{and}\quad  \overline{e}\in U(R_{m} ).
\]
\end{lemma}
\begin{proof} Let $\mu_a:=(a,m)$ and  $\mu_b:=(b,m)$.
This yields that
\[
\textstyle
\alpha _{b}=\frac{m}{\mu_b} =\frac{m}{\mu_a}   \frac{\mu_a}{\mu_b}=\alpha _{a}\frac{\mu_a}{\mu_b}
\]
and  $\overline{\alpha}_{b}=\overline{\alpha}_{a}\cdot \overline{\sigma}$ in which  $\sigma:=\frac{\mu_a}{\mu_b}$. Therefore,  $\overline{\sigma }+{\rm Ann}(\overline{\alpha}_{a})$ is the set of solutions of the equation $\overline{\alpha}_{b}=\overline{\alpha}_{a}\; \overline{x}$.
The ideal  ${\rm Ann}(\overline{\alpha}_{a})$ is generated by the image of $\alpha _{a}:=\frac{m}{\mu_a}\in R$ (see Lemma \ref{LL:5}), so
\[
\textstyle
\frac{m}{\left(\frac{m}{\mu_a} ,\; m\right)} =\frac{m}{\frac{m}{\mu_a} } =\mu_a\quad\text{and}\quad {\rm Ann}(\overline{\alpha}_{a})=\overline{\mu}_{a}R_{m}.
\]
Clearly  $\overline{\mu}_{a}=\overline{\mu}_{b}\; \overline{\sigma }$ and  $\overline{\sigma }| \overline{\mu}_{a}$.  This yields that $\overline{\sigma }$ is a generating solution  of the class $\overline{\sigma }+{\rm Ann}(\overline{\alpha}_{a})$. Consequently, \;
$\tifrac{\overline{a}}{\overline{b}} \;=\tifrac{\overline{\alpha} _{b} }{\overline{\alpha} _{a} } \cdot \overline{e}$\;  for some $\overline{e}\in U(R_{m})$.
\end{proof}


In order to simplify the notation, in the sequel of the paper we will omit the over line when referring to the elements of the ring $R_{m}$.

To a  permutation $\sigma=\left(\begin{matrix} {1} & {2} & {\ldots } & {n} \\ {i_{1} } & {i_{2} } & {\ldots } & {i_{n} } \end{matrix}\right)\in S_n$  we assign  the following two sets:
\[
\begin{split}
\mathfrak{I}_1(\sigma)&=\textstyle\{ \;
(p_{i}, q_{i})  \; \mid  \;  p_{i} >q_{i} \quad \text{and}   \; \binom{ p_{i}}{q_{i}} \; \text{is a column in }\; \sigma \};\\
\mathfrak{I}_2(\sigma)&=\textstyle\{ \;
(\alpha_{i}, \beta_{i})  \;\mid \;   \alpha_{i} \leq \beta_{i} \;\; \text{and}  \; \binom{ \alpha_{i}}{\beta_{i}} \; \text{is a column in }\;  \sigma \}.\\
\end{split}
\]

\begin{lemma}\label{LL:9}
If $\sigma\in S_n$, then
\begin{equation}\label{Eq:7}
\textstyle
\prod_{(p_i,q_i)\in \mathfrak{I}_1(\sigma)}
 \frac{p_{i} }{q_{i}}=
\prod_{(\alpha_i, \beta_i)\in \mathfrak{I}_2(\sigma)}
\frac{\beta _{i} }{\alpha _{i}}.
\end{equation}
Moreover, if $\Phi ={\rm diag}(\varphi _{1} ,\ldots ,\varphi _{n}) \in R_{m}^{n \times n}$ such that   $\varphi _{1} | \varphi _{2} |\cdots | \varphi _{n}\not=0$, then  each  set $\overline{M}_{ij}$  (see \eqref{Eq:5}) contains a generating solution  ${\psi_{ij}}$  (with $1 \leq j<i\leq  n$),  such that
\begin{equation}\label{Eq:8}
\prod_{(p_i,q_i)\in \mathfrak{I}_1(\sigma)}
{\psi}_{p_i,q_i}=
\prod_{(\alpha_i, \beta_i)\in \mathfrak{I}_2(\sigma)}
{\psi}_{\beta_i,\alpha_i}.
\end{equation}
\end{lemma}

\begin{proof}
Let us prove  \eqref{Eq:7}.
Let $\sigma_t=(i_1,i_2,\ldots,i_t)$ be  a cycle. The case $t=2$ is trivial. Now let  \eqref{Eq:7} holds for $\sigma_t$ where  $t>2$. Consider the cycle $\sigma_{t+1}=(i_1,i_2,\ldots,i_t, i_{t+1})$ of length $t+1$ which we obtain from $\sigma_t$ adding a new symbol $i_{t+1}$. Let us prove that the equation \eqref{Eq:7} holds for the new cycle $\sigma_{t+1}$, if we replace the pair $(i_t,i_1)$ by two pairs $(i_t,i_{t+1})$ and $(i_{t+1}, i_1)$ adding them  in appropriate places of the equation \eqref{Eq:7}. These three pairs are either lie in one set $\mathfrak{I}_{s_1}$ or two are in  $\mathfrak{I}_{s_2}$  and the one  is in $\mathfrak{I}_{s_3}$, where $s_1,s_2,s_3\in \{1,2\}$ and $s_2\not=s_3$.
Now considering these possible cases, it is easy to check that \eqref{Eq:7} holds for any cycle of finite length. Since each permutation is a product of disjoint cycles, \eqref{Eq:7} always holds.

The proof of \eqref{Eq:8} follows immediately  from \eqref{Eq:7},  Lemma \ref{LL:6}  and  Corollary \ref{Cor:2}.
\end{proof}

\begin{lemma}\label{LL:10}
If ${\psi_{ij}}$ is such  a generating solution of $\overline{M}_{ij}$ for all  $1 \leq j<i\leq  n$
for which  \eqref{Eq:8} holds, then
\[
\det \left[
\begin{smallmatrix} {h_{11} } & {h_{12} } & {\cdots} & {h_{1,\; n-1} } & {h_{1n} } \\
{\psi}_{21}h_{21} & {h_{22} } & {\cdots} & {h_{2,\; n-1} } & {h_{2n} } \\
 {\cdots} & {\cdots} & {\cdots} & {\cdots} & {\cdots} \\
{\psi}_{n1} h_{n1}  & {{\psi}_{n2} h_{n2} } & {\cdots} & {{\psi}_{n,n-1} h_{n,\; n-1} } & {h_{nn} } \end{smallmatrix}\right] =
\det \left[    \begin{smallmatrix} {h_{11} } & {{\psi}_{21} h_{12} } & {\cdots} & {{\psi}_{n-1,1} h_{1,\; n-1} } &
 {{\psi}_{n1} h_{1n} } \\
  {h_{21} } & {h_{22} } & {\cdots} & {{\psi}_{n-1,2} h_{2,\; n-1} } & {{\psi}_{n2} h_{2n} } \\
   {\cdots} & {\cdots} & {\cdots} & {\cdots} & {\cdots} \\ {h_{n1} } & {h_{n2} } & {\cdots} & {h_{n,\; n-1} } & {h_{nn} } \end{smallmatrix}
   \right],
   \]
in which  $h_{ij}\in R_m$.
\end{lemma}

\begin{proof}  Let us show that  both determinants consist of the same summands. Assign to each summand a permutation $\sigma=\left(\begin{smallmatrix} {1} & {2} & {\ldots } & {n} \\ {i_{1} } & {i_{2} } & {\ldots } & {i_{n} } \end{smallmatrix}\right)$. In the left hand side  determinant  these summands have the following form:
\[
(-1)^{{\rm sign}(\sigma)}{\psi}_{p_1,q_1}{\psi}_{p_2,q_2} \ldots {\psi}_{p_s,q_s} h_{p_{1}, q_{1} } \ldots h_{p_{s}, q_{s} } h_{\alpha _{1}, \beta _{1} } \ldots h_{\alpha _{t}, \beta_{t} },
\]
where $(p_1,q_1),\ldots, (p_s,q_s)\in \mathfrak{I}_{1}(\sigma)$,
$(\alpha _1,\beta_1),\ldots, (\alpha_t,\beta_t)\in \mathfrak{I}_{2}(\sigma)$, $s=|\mathfrak{I}_{1}(\sigma)|$ and  $t=|\mathfrak{I}_{2}(\sigma)|$.
The corresponding summands in the right hand side   determinant are
\[
(-1)^{{\rm sign}(\sigma)}{\psi}_{\beta_1,\alpha_1}{\psi}_{\beta_2,\alpha_2}  \ldots  {\psi}_{\beta_t,\alpha_t} h_{p_{1}, q_{1} } \ldots h_{p_{s}, q_{s} } h_{\alpha _{1}, \beta _{1} } \ldots h_{\alpha _{t}, \beta _{t} }.
\]
Finally, these terms  are the same  by \eqref{Eq:8} from Lemma \ref{LL:9}.\end{proof}

\begin{lemma}\label{LL:11}
Let
$\Phi ={\rm diag}(\varphi _{1} ,\ldots ,\varphi _{n}) \in R_{m}^{n \times n}$ such that   $\varphi _{1} | \varphi _{2} |\cdots | \varphi _{n}\not=0$.
If  $\alpha_{ij}$ and $\beta_{ij}$ are arbitrary elements from $\overline{M}_{ij}$ for all $1 \leq j<i\leq  n$, then
\[
 \left[
\begin{smallmatrix}{h_{11} } & {h_{12} } & {\cdots} & {h_{1,\; n-1} } & {h_{1n} } \\
\alpha_{21}h_{21} & {h_{22} } & {\cdots} & {h_{2, n-1} } & {h_{2n} } \\
 {\cdots} & {\cdots} & {\cdots} & {\cdots} & {\cdots} \\
\alpha_{n1} h_{n1}  & {\alpha_{n2} h_{n2} } & {\cdots} & {\alpha_{n,n-1} h_{n, n-1} } & {h_{nn} } \end{smallmatrix}
   \right]\Phi =
  \Phi  \left[
\begin{smallmatrix} {h_{11} } & {\beta_{21} h_{12} } & {\cdots} & {\beta_{n-1,1} h_{1, n-1} } &
 {\beta_{n1} h_{1n} } \\
  {h_{21} } & {h_{22} } & {\cdots} & {\beta_{n-1,2} h_{2, n-1} } & {\beta_{n2} h_{2n} } \\
   {\cdots} & {\cdots} & {\cdots} & {\cdots} & {\cdots} \\ {h_{n1} } & {h_{n2} } & {\cdots} & {h_{n, n-1} } & {h_{nn} } \end{smallmatrix}
   \right],
\]
in which $h_{ij}\in R_m$.
\end{lemma}

\begin{proof}
Each  element at  position $(i,j)$ where $i>j$ of the matrix from the left hand side  product  has  the  form $\varphi_{j}\alpha_{ij}h_{ij}$.
The corresponding element of the  matrix from the right  hand side  product has the form  $\varphi _{i}h_{ij}$.   Since  $\alpha_{ij}$ is a solution of the equation  $\varphi _{i}=\varphi _{j}x$,  we can replace  $\varphi _{i}$ by $\varphi _{j}\alpha_{ij}$.
This means that the elements of the matrices which are situated below the main diagonals coincide.

Let $i<j$.  Each element at  position  $(ij)$ of the matrix from the left hand side  product  has  the  form  $\varphi _{j}h_{ij}$.  The corresponding element of the  matrix from the right  hand side  product has the form $\varphi _{i}\beta_{ji}h_{ij}$, respectively.
As above, $\varphi _{j}=\varphi _{i}\beta_{ji}$. Taking into account that main diagonals of both  matrices are equal, we obtain the requested equality.
\end{proof}

\begin{proof}[Proof of Theorem \ref{Theorem:B}]
$\Rightarrow$. Let $H=[ p_{ij}] \in {\bf{G}}_{\Phi }$. By definition of the Zelisko group, there exists $S=[ s_{ij}]\in {\rm GL}_n(R_{m})$ such that $H\Phi =\Phi S$, so
\begin{equation}\label{Eq:9}
\varphi _{j} p_{ij} =\varphi _{i} s_{ij}, \qquad\qquad (1\leq i,j\leq n).
\end{equation}
Obviously, for $i\leq j$ no restrictions are imposed on the elements $p_ {ij}$.

Let  $i>j$.  Since  $\varphi _{j} \;|\; \varphi _{i}$,  we have   $\varphi _{i} =\varphi _{j}\cdot\tifracc{\varphi _{i} }{\varphi _{j}}$ by  Corollary \ref{Cor:1}. Now using  (\ref{Eq:9}) we obtain that
\[
\textstyle\varphi _{j} \left(p_{ij} -\tifracc{\varphi _{i} }{\varphi _{j} } s_{ij} \right)=0.
\]
So $q_{ij}:=p_{ij} -\tifracc{\varphi _{i} }{\varphi _{j} } s_{ij}  \in {\rm Ann}(\varphi _{j} )$.
This yields that  $p_{ij} \in \tifracc{\varphi _{i} }{\varphi _{j} }  s_{ij} +{\rm Ann}(\varphi _{j} )$ and ${\rm Ann}(\varphi _{j} )=\alpha _{\varphi _{j} } R_{m}$ by  Lemma \ref{LL:5}.
Furthermore,
$\tifracc{\varphi _{i} }{\varphi _{j} } =\tifraccc{\alpha _{\varphi _{j} } }{{\alpha _{\varphi _{i} } }}\;\cdot  e$ (see Lemma \ref{LL:8}),
in which   $e\in U(R_{m} )$. It follows that
\[\textstyle
\alpha _{\varphi _{j} }= \tifraccc{\alpha _{\varphi _{j} } }{{\alpha _{\varphi _{i} } }}\cdot \alpha _{\varphi _{i} }=(\tifracc{\varphi _{i} }{\varphi _{j} }\cdot e^{-1}) \alpha _{\varphi _i},
\]
by  Corollary \ref{Cor:1}, so $\tifracc{\varphi _{i} }{\varphi _{j} } | \alpha _{\varphi _{j} }$. Thus,  $\tifracc{\varphi _{i} }{\varphi _{j} }$ is a divisor of all elements of the ideal $\alpha _{\varphi _{j} } R_{m} ={\rm Ann}(\varphi _{j} )$.
This yields that  $\tifracc{\varphi _{i} }{\varphi _{j} }| q_{ij}$ and $q_{ij} =\tifracc{\varphi _{i} }{\varphi _{j} }\cdot l_{ij}$ for some $l_{ij}\in R_m$.
Hence,
\[
\textstyle
p_{ij} = \tifracc{\varphi _{i} }{\varphi _{j} } s_{ij}+q_{ij}=
 \tifracc{\varphi _{i} }{\varphi _{j} } s_{ij}+\tifracc{\varphi _{i} }{\varphi _{j} } l_{ij} =
 \tifracc{\varphi _{i} }{\varphi _{j} } (s_{ij}+ l_{ij})
 = \tifracc{\varphi _{i} }{\varphi _{j} } h_{ij},
\]
in which $h_{ij}:=s_{ij}+ l_{ij}$.
Therefore,  the matrix $H=[ p_{ij}]$ has the form \eqref{Eq:2}.

$\Leftarrow$. The sequential product of  generating solutions of the sets $\overline{M}_{i,i-1}$, $\overline{M}_{i-1,i-2}, \ldots, \overline{M}_{j+1,j}$ (see \eqref{Eq:5}) is denoted  by
\[
\textstyle \psi_{ij}:={\tifracccc{\varphi _{i} }{\varphi _{i-1} }}\cdot { \tifracccc{\varphi _{i-1} }{\varphi _{i-2} }} \cdots{\tifracccc{\varphi _{j+2} }{\varphi _{j+1} }}\cdot { \tifracccc{\varphi _{j+1} }{\varphi _{j} }}, \qquad\qquad  (1\leq j<i\leq n).
 \]
Each   $\psi_{ij} $ is  the generating solution of $\overline{M}_{ij}$ by Lemma \ref{LL:7}.  Moreover,
$ {\tifracc{\varphi _{i} }{\varphi _{j} }}= \psi_{ij}e_{ij}$ in which   $e_{ij} \in U(R_m)$.
Hence,  we have a presentation
\[
H=  \left[
\begin{smallmatrix} {h_{11} } & {h_{12} }  & {h_{13} } & {\cdots} & {h_{1, n-2} } & {h_{1, n-1} } & {h_{1n} } \\
{\psi}_{21}e_{21}h_{21} & {h_{22}}  & {h_{23} } & {\cdots} & {h_{2, n-2} } & {h_{2, n-1} } & {h_{2n} } \\
{\psi}_{31}e_{31} h_{31} & {\psi}_{32}e_{32}{h_{32} }  & {h_{33} } & {\cdots} & {h_{3, n-2} }  & {h_{3, n-1} } & {h_{3n} } \\
 {\cdots} & {\cdots} & {\cdots} & {\cdots} & {\cdots} & {\cdots}  & {\cdots}\\
{\psi}_{n1}e_{n1} h_{n1}  & {{\psi}_{n2}e_{n2} h_{n2} }& {\cdots}  & {\cdots} & {\psi}_{n,n-2}e_{n,n-2}{h_{n, n-2} } & {\psi}_{n,n-1}e_{n,n-1} h_{n,n-1}  & {h_{nn} } \end{smallmatrix}
   \right].
\]
According to Lemma \ref{LL:10}, the determinants of the matrix $H$  and the matrix
 \[  H_1=  \left[
\begin{smallmatrix}
{h_{11} } & {\psi}_{21} {h_{12} }  & {\psi}_{31}{h_{13} } & {\cdots} & {\psi}_{n-2,1}{h_{1, n-2} } & {\psi}_{n-1,1}{h_{1, n-1} } & {\psi}_{n1}{h_{1n} } \\
e_{21}h_{21} & {h_{22} }  & {\psi}_{32}{h_{23} } & {\cdots} &  {\psi}_{n-2, 2}{h_{2, n-2} } &  {\psi}_{n-1, 2}{h_{2, n-1} } & {\psi}_{n2}{h_{2n} } \\
e_{31} h_{31} & e_{32}{h_{32} }  & {h_{33} } & {\cdots} &  {\psi}_{n-2, 3}{h_{3, n-2} }  & {\psi}_{n-1, 3}{h_{3, n-1} } & {\psi}_{n3}{h_{3n} } \\
 {\cdots} & {\cdots} & {\cdots} & {\cdots} & {\cdots} & {\cdots}  & {\cdots}\\
e_{n1} h_{n1}  & e_{n2} h_{n2} & {\cdots}  & {\cdots} & e_{n,n-2}{h_{n, n-2} } &  e_{n,n-1}h_{n,n-1}  & {h_{nn} } \end{smallmatrix}
   \right]
\]
coincide, so  $H_1$ is invertible. Now, using Lemma \ref{LL:11} gives     $H\Phi=\Phi H_1$, so   $H \in {\bf{G}}_{\Phi }$.\end{proof}

We propose  the following.
\bigskip

\noindent
{\bf Problem.}
Describe that rings $R$ in which g.c.d.  of all solutions of a solvable  linear equation   $b=ax$ ($a,b\in R$)  in $R$  is again a solution of the same  linear equation.

Note that, for rings  $M_n(R)$ over  elementary divisor domains   $R$ a positive solution to this problem   was done in \cite{Shch_New}.

\section{Acknowledgement}
\noindent
Authors would like to express their gratitude to the referee for valuable remarks. The work was supported by the UAEU UPAR [grant number G00002160].

\bibliographystyle{abbrv}

\end{document}